\documentclass[12pt]{amsart}

\newtheorem{theorem}{Theorem}[section]

\newtheorem{corollary}[theorem]{Corollary}
\theoremstyle{definition}

\theoremstyle{remark}
\newtheorem{remark}[theorem]{Remark}
\numberwithin{equation}{section}
\newfont{\kh}{msbm10 at 12pt}
\newcommand{\R}{\mbox{\kh R}}
\newcommand{\C}{\mbox{\kh C}}

\begin{document}
\title{A Fixed Point Approach to Stability of a Quadratic Equation}
\author{M. Mirzavaziri and M. S. Moslehian}
\address{Madjid Mirzavaziri\newline Department of Mathematics, Ferdowsi University, P. O. Box 1159, Mashhad 91775, Iran}
\email{mirzavaziri@math.um.ac.ir}
\address{Mohammad Sal Moslehian\newline Department of Mathematics, Ferdowsi University, P. O. Box 1159, Mashhad 91775, Iran}
\email{moslehian@ferdowsi.um.ac.ir} 
\begin{abstract}
Using the fixed point alternative theorem we establish the
orthogonal stability of quadratic functional equation of Pexider
type $f(x+y)+g(x-y)=h(x)+k(y)$, where $f, g, h, k$ are mappings
from a symmetric orthogonality space to a Banach space, by
orthogonal additive mappings under a necessary and sufficient
condition on $f$.
\end{abstract}
\maketitle

\section{Introduction.}
Suppose that ${\mathcal X}$ is a real vector space with $\dim
{\mathcal X}\geq 2$ and $\perp$ is a binary
relation on ${\mathcal X}$ with the following properties:\\
(O1) \textit{totality of $\perp$ for zero}: $x\perp 0, 0\perp x$ for all $x\in {\mathcal X}$;\\
(O2) \textit{independence}: if $x,y\in {\mathcal X}-\{0\}, x\perp y$, then $x,y$ are linearly independent;\\
(O3) \textit{homogeneity}: if $x,y\in {\mathcal X}, x\perp y$, then $\alpha x\perp \beta y$ for all $\alpha, \beta$ in the real line $\R$;\\
(O4) \textit{the Thalesian property}: Let $P$ be a $2$-dimensional
subspace of ${\mathcal X}$. If $x\in P$ and $\lambda$ in the
nonnegative real numbers $\R_+$, then there exists $y_0\in P$
such that $x\perp y_0$ and $x+y_0\perp \lambda x-y_0$.

Then the pair $({\mathcal X},\perp)$ is called an
\textit{orthogonality space}; cf. \cite{RAT}. By an
\textit{orthogonality normed space} we mean an orthogonality
space equipped with a norm. Some examples of special interest are

(i) The trivial orthogonality on a vector space ${\mathcal X}$
defined by (O1), and for non-zero elements $x,y\in {\mathcal X}$,
$x\perp y$ if and only if $x,y$ are linearly independent.

(ii) The ordinary orthogonality on an inner product space
$({\mathcal X}, \langle.,.\rangle)$ given by $x\perp y$ if and
only if $\langle x,y\rangle=0$.

(iii) The Birkhoff-James orthogonality on a normed space
$({\mathcal X},\|.\|)$ defined by $x\perp y$ if and only if
$\|x+\lambda y\|\geq \|x\|$ for all $\lambda\in \R$; cf.
\cite{JAM}.

The relation $\perp$ is called \textit{symmetric} if $x\perp y$
implies that $y\perp x$ for all $x,y\in {\mathcal X}$. Clearly
examples (i) and (ii) are symmetric but example (iii) is not. It
is remarkable to note, however, that a real normed space of
dimension greater than or equal to $3$ is an inner product space
if and only if the Birkhoff-James orthogonality is symmetric; see
\cite{AMI}.

Let ${\mathcal X}$ be a vector space (an orthogonality space) and
$({\mathcal G},+)$ be an abelian group. A mapping $f:{\mathcal
X}\to {\mathcal G}$ is called \textit{(orthogonally) additive} if
it satisfies the so-called \textit{(orthogonal) additive
functional equation} $f(x+y)=f(x)+f(y)$ for all $x,y\in {\mathcal
X}$ (with $x\perp y$). A mapping $f:{\mathcal X}\to {\mathcal G}$
is said to be \textit{(orthogonally) quadratic} if it satisfies
the so-called \textit{(orthogonally) Jordan-von Neumann quadratic
functional equation} $f(x+y)+f(x-y)=2f(x)+2f(y)$ for all $x,y\in
{\mathcal X}$ (with $x\perp y$).

The problem of  ``stability of functional equations'' is that
``when the solutions of an equation differing slightly from a
given one must be close to an exact solution of the given
equation?''. In 1941, S. M. Ulam \cite{ULA} posed the first
question on the subject concerning the stability of group
homomorphisms. In 1941, D. H. Hyers \cite{HYE} gave a partial
solution of Ulam's problem in the context of Banach spaces. In
1978, Th. M. Rassias \cite{RAS1} generalized the theorem of Hyers
to an unbounded situation. The result of Rassias has provided a
lot of influence in the development of what we now call {\it
Hyers--Ulam--Rassias stability} of functional equations.
Following Hyers and Rassias approaches, during the last decades,
the stability problem for several functional equations have been
extensively investigated by many mathematicians; cf.
\cite{H-I-R}. Nowadays, there may be found several applications
in actuarial and financial mathematics, sociology, psychology,
and pure mathematics \cite{A-D}.

The first author who treated the stability of the quadratic
equation was F. Skof \cite{SKO}. P. W. Cholewa \cite{CHO}
extended Skof's theorem to abelian groups. Skof's result was also
generalized by S. Czerwik \cite{CZE1} in the spirit of
Hyers--Ulam--Rassias. S. M. Jung \cite{JUN1, JUN2} investigated
the stability of the quadratic equation. K.W. Jun and Y. H. Lee
\cite{J-L} proved the stability of quadratic equation of Pexider
type. The stability problem of the quadratic equation has been
extensively investigated by some mathematicians; cf. \cite{CZE2,
CZE3, RAS2}.

The orthogonal quadratic equation $$f(x+y)+f(x-y)=2f(x)+2f(y),~
x\perp y$$ was first investigated by F. Vajzovi\' c \cite{VAJ}
when ${\mathcal X}$ is a Hilbert space, ${\mathcal G}$ is the
scalar field, $f$ is continuous and $\perp$ means the Hilbert
space orthogonality. H. Drljevi\' c \cite{DRL2} proved the
following stability result:

{\small Let ${\mathcal H}$ be a complex Hilbert space of
dimension $\geq 3$, and $A:{\mathcal H}\to {\mathcal H}$ a bounded
self-adjoint linear operator with $\dim A({\mathcal H})\geq 2$,
and let the real numbers $\theta\geq 0$ and $p\in [0,2)$ be given.
Suppose that $f:{\mathcal H}\to\C$ is continuous and satisfies the
inequality
\begin{eqnarray*}
|f(x+y)+f(x-y)-2f(x)-2f(y)|\leq \theta[|\langle
x,x\rangle|^{p/2}+|\langle y,y\rangle|^{p/2}],
\end{eqnarray*}
whenever $\langle Ax,y\rangle=0$. Then the limit
$T(x)=\lim_{n\to\infty}\frac{f(2^nx)}{4^n}$ exists for each $x\in
{\mathcal H}$ and the functional $T$ is continuous and satisfies
$T(x+y)+T(x-y)=2T(x)+2T(y)$ whenever $\langle Ax,y\rangle=0$.
Moreover, there exists a real number $\varepsilon>0$ such that
\begin{eqnarray*}
|f(x)-T(x)|\leq\varepsilon|\langle Ax,x\rangle|^{p/2},
\end{eqnarray*} for all $x\in {\mathcal H}$.}

Later H. Drljevi\' c \cite{DRL1}, M. Fochi \cite{FOC} and G.
Szab\' o \cite{SZA} obtained more results on the subject.

One of the significant conditional equations is the so-called
\textit{orthogonally quadratic functional equation of Pexider
type} $f(x+y)+g(x-y)=h(x)+k(y),~~ x\perp y$. Recently, the second
author investigated this equation with ``$g=f$''. Using some ideas
from \cite{G-S, J-S, MOS1, MOS2, J-C, C-R}, we aim to use the
alternative of fixed point theorem to establish the stability of
this equation in the spirit of Hyers--Ulam under certain
conditions. The first systematic study of fixed point theorems in
nonlinear analysis is due to G. Isac and Th. M. Rassias; cf.
\cite{I-R}.

\section{Main Results.}
We start our work with a known fixed point theorem which will be
needed later:

\begin{theorem} (The alternative of fixed point) Suppose $({\mathcal E},d)$ be a complete generalized metric space and $J:{\mathcal E}\to {\mathcal E}$ be a strictly contractive mapping with the Lipschitz constant L. Then, for each given element $x\in {\mathcal E}$,either\\
(A1) $d(J^nx, J^{n+1}x)=\infty$\\
for all $n\geq 0$, or\\
(A2) There exists a natural number $n_0$ such that:\\
(A20) $d(J^nx, J^{n+1}x)<\infty$, for all $n\geq n_0$;\\
(A21) The sequence $\{J^nx\}$ is convergent to a fixed point $y^*$ of $J$;\\
(A22) $y^*$ is the unique fixed point of $J$ in the set $Y=\{y\in {\mathcal E}: d(J^{n_0}x,y)<\infty\}$;\\
(A23) $d(y,y^*)<\frac{1}{1-L}d(y,Jy)$ for all $y\in
Y$.\end{theorem}

Suppose that ${\mathcal X}$ denotes an orthogonality real space
and ${\mathcal Y}$ denotes a Banach space. Consider the set
${\mathcal E}:=\{\varphi:{\mathcal X}\to {\mathcal Y}:
\varphi(0)=0\}$ and introduce a generalized metric on ${\mathcal
E}$~ by
\begin{eqnarray*}
d(\varphi,\psi)=\inf\{c\in (0,\infty): \|\varphi(x)-\psi(x)\|\leq
c, \forall x\in {\mathcal X}\}.
\end{eqnarray*}
It is easy to see that $({\mathcal E},d)$ is complete. Given a
number $0\leq \lambda<1$, define the following mapping
$J_\lambda:{\mathcal E}\to {\mathcal E}$ by $(J_\lambda
\varphi)(x):=\lambda \varphi(2x)$. For arbitrary elements
$\varphi, \psi\in {\mathcal E}$ we have
\begin{eqnarray*}
d(\varphi,\psi)<c &\Rightarrow& \|\varphi(x)-\psi(x)\|\leq c,~~~x\in {\mathcal X}\\
&\Rightarrow& \|\lambda \varphi(2x)-\lambda \psi(2x)\|\leq\lambda c,~~~x\in {\mathcal X}\\
&\Rightarrow& d(J_\lambda \varphi,J_\lambda \psi)\leq\lambda c.
\end{eqnarray*}

Therefore
\begin{eqnarray*}
d(J_\lambda \varphi,J_\lambda \psi)\leq\lambda
d(\varphi,\psi),~~~\varphi,\psi\in {\mathcal E}.
\end{eqnarray*}
Hence $J_\lambda$ is a strictly contractive mapping on ${\mathcal
E}$ with the Lipschitz constant $\lambda$ and we can use the fixed
point alternative theorem.

We are just ready to prove the orthogonal stability of the
Pexiderized equation $f(x+y)+g(x+y)=h(x)+k(y)$ where $f, g, h, k$
are mappings from ${\mathcal X}$ to ${\mathcal Y}$ under certain
condition.

We use the notation $\varphi(x)\unlhd~\varepsilon$ in the sense
that there exists a number $a$ such that $\varphi(x)\leq
a\varepsilon$ for all $x$ in the domain of $\varphi$.

\begin{theorem} Suppose that ${\mathcal X}$ is a real orthogonality space with
a symmetric orthogonal relation $\perp$ and ${\mathcal Y}$ is a
Banach space. Let the mappings $f, g, h, k:{\mathcal X}\to
{\mathcal Y}$ satisfy the following inequalities
\begin{eqnarray}
\|f(x+y)+g(x-y)-h(x)-k(y)\|\leq\varepsilon,
\end{eqnarray}
for all $x, y\in {\mathcal X}$ with $x\perp y$. Then there exists
an orthogonally additive mapping $T$ such that
\begin{eqnarray*}
\|f(x)-T(x)\|\unlhd~\varepsilon
\end{eqnarray*}
if and only if
\begin{eqnarray*}
\|f(2x)-f(-2x)-4f(x)-4f(-x)\|\unlhd~\varepsilon.
\end{eqnarray*}
Indeed, if
\begin{eqnarray}
\|f(2x)-f(-2x)-4f(x)-4f(-x)\|\leq\varepsilon.
\end{eqnarray}
holds for all $x\in {\mathcal X}$, then there exist orthogonally
additive mappings $T, T', T'':{\mathcal X}\to {\mathcal Y}$ such
that
\begin{eqnarray*}
\|f(x)-f(0)-T(x)\|\leq \frac{140}{3}\varepsilon,
\end{eqnarray*}
\begin{eqnarray*}
\|g(x)-g(0)-T'(x)\|\leq\frac{98}{3}\varepsilon,
\end{eqnarray*}
\begin{eqnarray*}
\|h(x)+k(x)-h(0)-k(0)-T''(x)\|\leq\frac{256}{3}\varepsilon,
\end{eqnarray*}
for all $x\in {\mathcal X}$.
\end{theorem}
\begin{proof}
Suppose that $(2.2)$ holds. Define $F(x)=f(x)-f(0),
G(x)=g(x)-g(0), H(x)=h(x)-h(0), K(x)=k(x)-k(0)$. Then
$F(0)=G(0)=H(0)=K(0)=0$. Set $L(x)=\frac{H(x)+K(x)}{2}$.

Use (O1) and put $x=y=0$ in $(2.1)$ and subtract the argument of
the norm of the resulting inequality from that of inequality
$(2.1)$ to get
\begin{eqnarray}
\|F(x+y)+G(x-y)-H(x)-K(y)\|\leq 2\varepsilon.
\end{eqnarray}
Let $\rho^e(x)=\frac{\rho(x)+\rho(-x)}{2}$ and
$\rho^o(x)=\frac{\rho(x)-\rho(-x)}{2}$ denote the even and odd
parts of a given function $\rho$, respectively.

If $x\perp y$ then, by (O3), $-x\perp -y$. Hence we can replace $x$
by $-x$ and $y$ by $-y$ in $(2.3)$ to obtain
\begin{eqnarray}
\|F(-x-y)+G(-x+y)-H(-x)-K(-y)\|\leq 2\varepsilon.
\end{eqnarray}
By virtue of triangle inequality and $(2.3)$ and $(2.4)$ we have
\begin{eqnarray}
\|F^o(x+y)+G^o(x-y)-H^o(x)-K^o(y)\|\leq 2\varepsilon,
\end{eqnarray}
\begin{eqnarray}
\|F^e(x+y)+G^e(x-y)-H^e(x)-K^e(y)\|\leq 2\varepsilon,
\end{eqnarray}
for all $x,y\in {\mathcal X}$.\\
{\bf Step (I). Approximating $F^o$}

Let $x\perp y$. Then $y\perp x$, and by (2.5)
\begin{eqnarray}
\|F^o(x+y)-G^o(x-y)-H^o(y)-K^o(x)\|\leq 2\varepsilon.
\end{eqnarray}
It follows from $(2.5)$ and $(2.7)$ that
\begin{eqnarray}
&&\|2F^o(x+y)-H^o(x)-K^o(x)-H^o(y)-K^o(y)\|\nonumber\\&\leq&\|F^o(x+y)+G^o(x-y)-H^o(x)-K^o(y)\|\nonumber\\
&&+\|F^o(x+y)-G^o(x-y)-H^o(y)-K^o(x)\|\nonumber\\
&\leq& 4\varepsilon.
\end{eqnarray}
for all  $x,y\in {\mathcal X}$ with $x\perp y$. In particular, for
arbitrary $x$ and $y=0$ we get
\begin{eqnarray}
\|2F^o(x)-H^o(x)-K^o(x))\|&\leq& 4\varepsilon.
\end{eqnarray}
By $(2.8)$ and $(2.9)$, we have
\begin{eqnarray}
\|F^o(x+y)-F^o(x)-F^o(y)\|&\leq&\frac{1}{2}\|2F^o(x+y)-H^o(x)-K^o(x)-H^o(y)-K^o(y)\|\nonumber\\
&&+\frac{1}{2}\|2F^o(x)-H^o(x)-K^o(x))\|\nonumber\\
&&+\frac{1}{2}\|2F^o(y)-H^o(y)-K^o(y))\|\nonumber\\
&\leq& 6\varepsilon
\end{eqnarray}
for all  $x,y\in {\mathcal X}$ with $x\perp y$.

Given $x\in {\mathcal X}$, by (O4), there exists $y_0\in
{\mathcal X}$ such that $x\perp y_0$ and $x+y_0\perp x-y_0$.
Replacing $x$ and $y$ by $x+y_0$ and $x-y_0$ in $(2.10)$, we have
\begin{eqnarray}
\|F^o(2x)-F^o(x+y_0)-F^o(x-y_0)\|\leq 6\varepsilon.
\end{eqnarray}
Since $x\perp y_0$ and $x\perp -y_0$, it follows from $(2.10)$ that
\begin{eqnarray}
\|F^o(x+y_0)-F^o(x)-F^o(y_0)\|\leq 6\varepsilon,
\end{eqnarray}
and
\begin{eqnarray}
\|F^o(x-y_0)-F^o(x)+F^o(y_0)\|\leq 6\varepsilon.
\end{eqnarray}
By $(2.11), (2.12)$ and $(2.13)$,
\begin{eqnarray*}
\|\frac{1}{2}F^o(2x)-F^o(x)\|&\leq&\frac{1}{2}\|F^o(2x)-F^o(x+y_0)-F^o(x-y_0)\|\\
&&+\frac{1}{2}\|F^o(x+y_0)-F^o(x)-F^o(y_0)\|\\
&&+\frac{1}{2}\|F^o(x-y_0)-F^o(x)+F^o(y_0)\|\\
&\leq& 9\varepsilon.
\end{eqnarray*}
Hence $d(F^o,J_{1/2}F^o)\leq 9\varepsilon<\infty$. Using the fixed
point alternative we conclude the existence of a mapping
$R:{\mathcal X}\to {\mathcal Y}$ such that $R$ is a fixed point of
$J_{1/2}$ that is $R(2x)=2R(x)$ for all $x\in {\mathcal X}$. Since
$\displaystyle{\lim_{n\to\infty}}d(J_{1/2}^nF^o,R)=0$ we easily
deduce that
$\displaystyle{\lim_{n\to\infty}}\frac{F^o(2^nx)}{2^n}=R(x)$ for
all $x\in {\mathcal X}$.

Indeed, the mapping $R$ is the unique fixed point of $J_{1/2}$ in
the set $Y=\{\varphi\in {\mathcal E}: d(F^o,\varphi)<\infty\}$.
Hence $R$ is the unique fixed point of $J_{1/2}$ such that
$\|F^o(x)-R(x)\|\leq K$ for some $K>0$ and for all $x\in
{\mathcal X}$. Again, by applying the fixed point alternative
theorem we obtain
\begin{eqnarray*}
d(F^o,R)\leq 2d(F^o,J_{1/2}F^o)\leq 18\varepsilon.
\end{eqnarray*}
Thus
\begin{eqnarray}
\|F^o(x)-R(x)\|\leq 18\varepsilon,
\end{eqnarray}
for all $x\in {\mathcal X}$. Let $x\perp y$ and $n$ be a positive
integer. Then $2^nx\perp 2^ny$ and so we can replace $x$ and $y$
in $(2.10)$ by $2^nx$ and $2^ny$, respectively. Dividing the both
sides by $2^n$ and letting $n$ tend to $\infty$ we infer that
$R(x+y)=R(x)+R(y)$.
Hence $R$ is orthogonally additive.\\
{\bf Step (II). Approximating $G^o$}

Let $x\perp y$. Then $x\perp -y$ and $(2.5)$ yields the following.
\begin{eqnarray*}
\|G^o(x+y)+F^o(x-y)-H^o(x)-(-K)^o(y)\|\leq 2\varepsilon.
\end{eqnarray*}
Using the same argument as in Step (I), we conclude that there
exists a unique orthogonally additive mapping $R':{\mathcal X}\to
{\mathcal Y}$ such that
\begin{eqnarray}
\|G^o(x)-R'(x)\|\leq 18\varepsilon.
\end{eqnarray}
{\bf Step (III). Approximating $L^o$}

Using $2.9$ we get
\begin{eqnarray}
\|F^o(x)-L^o(x)\|\leq 2\varepsilon,
\end{eqnarray}
so ,by $(2.14)$,
\begin{eqnarray}
\|L^o(x)-R(x)\|\leq\|F^o(x)-L^o(x)\|+\|F^o(x)-R(x)\|\leq
2\varepsilon+ 18\varepsilon= 20\varepsilon.
\end{eqnarray}
{\bf Step (IV). Approximating $G^e$}

Now we use inequality $(2.6)$ concerning the even parts. Let
$x\perp y$. Then $y\perp x$ and by $(2.6)$ we get
\begin{eqnarray}
\|F^e(x+y)+G^e(x-y)-H^e(y)-K^e(x)\|\leq 2\varepsilon.
\end{eqnarray}
By $(2.6)$ and $(2.18)$ we infer that
\begin{eqnarray}
\|F^e(x+y)+G^e(x-y)-L^e(x)-L^e(y)\|\leq 2\varepsilon,
\end{eqnarray}
for all $x,y\in {\mathcal X}$ with $x\perp y$. In particular, it
follows from $x\perp 0$ that
\begin{eqnarray}
\|F^e(x)+G^e(x)-L^e(x)\|\leq 2\varepsilon,
\end{eqnarray}
for all $x\in {\mathcal X}$. Applying $(2.19)$ and $(2.20)$ we get
\begin{eqnarray}
&&\|\big(F^e(x+y)-F^e(x)-F^e(y)\big)+\big(G^e(x-y)-G^e(x)-G^e(y)\big)\|\nonumber\\
&\leq&\|F^e(x+y)+G^e(x-y)-L^e(x)-L^e(y)\|\nonumber\\
&&+\|F^e(x)+G^e(x)-L^e(x)\|+\|F^e(y)+G^e(y)-L^e(y)\|\nonumber\\
&\leq& 6\varepsilon,
\end{eqnarray}
for all $x,y\in {\mathcal X}$ with $x\perp y$.

Given $x\in {\mathcal X}$, by (O4), there exists $y_0\in
{\mathcal X}$ such that $x\perp y_0$ and $x+y_0\perp x-y_0$.
Hence, by (O3), $x\perp -y_0$ $x+y_0\perp y_0-x$ and so, by
repeatedly applying $(2.21)$, we get
\begin{eqnarray}
&&\|\big(F^e(x+y_0)-F^e(x)-F^e(y_0)\big)+\big(G^e(x-y_0)-G^e(x)-G^e(y_0)\big)\|\nonumber\\
&\leq&6\varepsilon,
\end{eqnarray}
\begin{eqnarray}
&&\|\big(F^e(x-y_0)-F^e(x)-F^e(y_0)\big)+\big(G^e(x+y_0)-G^e(x)-G^e(y_0)\big)\|\nonumber\\
&\leq&6\varepsilon,
\end{eqnarray}
\begin{eqnarray}
\|\big(F^e(2y_0)-F^e(x+y_0)-F^e(x-y_0)\big)+\big(G^e(2x)-G^e(x+y_0)-G^e(x-y_0)\big)\|\nonumber\\
&\leq&6\varepsilon.
\end{eqnarray}
By (O3), $\frac{x+y_0}{2}\perp\pm\frac{x-y_0}{2}$ and so by using
$(2.21)$, we obtain
\begin{eqnarray}
\|\big(F^e(x)-F^e(\frac{x+y_0}{2})-F^e(\frac{x-y_0}{2})\big)+\big(G^e(y_0)-G^e(\frac{x+y_0}{2})-G^e(\frac{x-y_0}{2})\big)\|\nonumber\\
&\leq& 6\varepsilon,
\end{eqnarray}
and
\begin{eqnarray}
\|\big(F^e(y_0)-F^e(\frac{x+y_0}{2})-F^e(\frac{x-y_0}{2})\big)+\big(G^e(x)-G^e(\frac{x+y_0}{2})-G^e(\frac{x-y_0}{2})\big)\|\nonumber\\
&\leq& 6\varepsilon.
\end{eqnarray}
It follows from $(2.25)$ and $(2.26)$ we infer that
\begin{eqnarray}
\|(F^e(x)-F^e(y_0))-(G^e(x)-G^e(y_0))\|\leq 12\varepsilon.
\end{eqnarray}
Using the triangle inequality, we infer from $(2.22), (2.23), (2.24)$ and $(2.27)$ that
\begin{eqnarray}
\|(F^e(2y_0)-4F^e(y_0))+(G^e(2x)-4G^e(x)\|\leq 42\varepsilon.
\end{eqnarray}
So far, we do not use $(2.2)$. Now we may apply $(2.2)$ and
$(2.28)$ to get
\begin{eqnarray*}
\|\frac{1}{4}G^e(2x)-G^e(x)\|&\leq&\frac{1}{4}\|F^e(2y_0)-4F^e(y_0)\|+ \frac{42}{4}\varepsilon\\
&\leq&\frac{\varepsilon}{2}+\frac{21\varepsilon}{2}=11\varepsilon.
\end{eqnarray*}
Therefore $d(G^e,J_{1/4}G^e)\leq 11\varepsilon<\infty$. Using the
fixed point alternative we conclude the existence of a mapping
$S':{\mathcal X}\to {\mathcal Y}$ such that $S'$ is a fixed point
of $J_{1/4}$ that is $S'(2x)=4S'(x)$ for all $x\in {\mathcal X}$.
Since $\displaystyle{\lim_{n\to\infty}}d(J_{1/4}^nG^e,S')=0$ we
easily deduce that
$\displaystyle{\lim_{n\to\infty}}\frac{G^e(2^nx)}{2^{2n}}=S'(x)$
for all $x\in {\mathcal X}$.

Indeed, the mapping $S'$ is the unique fixed point of $J_{1/4}$ in
the set $Y=\{\psi\in {\mathcal E}: d(G^e,\psi)<\infty\}$. Hence
$S'$ is the unique fixed point of $J_{1/4}$ such that
$\|G^e(x)-S'(x)\|\leq K$ for some $K>0$ and for all $x\in
{\mathcal X}$. Again, by applying the fixed point alternative
theorem we obtain
\begin{eqnarray*}
d(G^e,S')\leq\frac{4}{3}d(G^e,J_{1/4}G^e)\leq\frac{44}{3}\varepsilon.
\end{eqnarray*}
Thus
\begin{eqnarray}
\|G^e(x)-S'(x)\|\leq \frac{44}{3}\varepsilon.
\end{eqnarray}
Let $x\perp y$ and $n$ be a positive integer. Then $2^nx\perp
2^ny$ and so we can replace $x$ and $y$ in $(2.10)$ by $2^nx$ and
$2^ny$, respectively. Dividing the both sides by $2^{2n}$ and
taking the limit as $n\to\infty$ we infer that
$S'(x+y)=S'(x)+S'(y)$. Hence $S'$ is orthogonally additive.

{\bf Step (V). Approximating $F^e$}

Let $x\perp y$. Then $x\perp -y$ and $(2.6)$ yields the following.
\begin{eqnarray*}
\|G^e(x+y)+F^e(x-y)-H^e(x)-K^e(y)\|\leq 2\varepsilon.
\end{eqnarray*}
By $(2.28)$, we have
\begin{eqnarray}
\|G^e(2x)-4G^e(x)\|&\leq&\|F^e(2y_0)-4F^e(y_0)\|+
42\varepsilon\leq 44\varepsilon.
\end{eqnarray}
Using the same argument as in Step (IV) and noting to $(2.30)$,
we conclude the existence of a unique orthogonally additive
mapping $S:{\mathcal X}\to {\mathcal Y}$ such that
$S(x)=\lim_{n\to\infty}\frac{F^e(2^nx)}{2^{2n}}$ and
\begin{eqnarray}
\|F^e(x)-S(x)\|\leq \frac{86}{3}\varepsilon.
\end{eqnarray}

{\bf Step (VI). Approximating $L^e$}

Inequalities $(2.20), (2.29)$ and $(2.31)$ yield the following.
\begin{eqnarray}
\|L^e(x)-S(x)-S'(x)\|&\leq&
\|F^e(x)+G^e(x)-L^e(x)\|+\|F^e(x)-S(x)\|-\|G^e(x)-S'(x)\|\nonumber\\
&\leq&
2\varepsilon+\frac{86}{3}+\frac{44}{3}\varepsilon\varepsilon\nonumber\\
&=&\frac{136}{3}\varepsilon.
\end{eqnarray}

{\bf Step (VII). Approximating $f, g, h+k$}

Put $T(x)=R(x)+S(x), T'(x)=R'(x)+S'(x)$ and
$T''(x)=2R(x)+2S(x)+2S'(x)$. Then $T, T'$ and $T''$ are
orthogonally additive and $(2.14), (2.15), (2.17), (2.29),
(2.31)$ and $(2.32)$ yield the following inequalities for each
$x\in {\mathcal X}$:
\begin{eqnarray*}
\|f(x)-f(0)-T(x)\|\leq \|F^o(x)-R(x)\|+\|F^e(x)-S(x)\|\leq
18\varepsilon+\frac{86}{3}\varepsilon=\frac{140}{3}\varepsilon,
\end{eqnarray*}
\begin{eqnarray*}
\|g(x)-g(0)-T'(x)\|\leq \|G^o(x)-R'(x)\|+\|G^e(x)-S'(x)\|\leq
18\varepsilon+\frac{44}{3}\varepsilon=\frac{98}{3}\varepsilon,
\end{eqnarray*}
\begin{eqnarray*}
\|h(x)+k(x)-h(0)-k(0)-T''(x)\|&\leq& 2\|L^o(x)-R(x)\|+2\|L^e(x)-S(x)-S'(x)\|\\
&\leq& 40\varepsilon+\frac{136}{3}\varepsilon\\
&=&\frac{256}{3}\varepsilon.
\end{eqnarray*}

{\bf Step (VIII). Necessity}

Let $T$ be an orthogonally additive mapping such that
$\|f(x)-T(x)\|\unlhd~\varepsilon$. Then
$\|f^e(x)-T^e(x)\|\unlhd~\varepsilon$. Note that $T^e$ is an
orthogonally additive mapping.

Let $x\in {\mathcal X}$. Using (O4), there exists a vector
$y_0\in {\mathcal X}$ such that $x\perp y_0$ and $x+y_0\perp
x-y_0$. Then, by (O3), $\frac{x}{2}\perp\frac{y_0}{2},
\frac{x+y_0}{2}\perp\frac{x-y_0}{2}$ and $x+y_0\perp y_0-x$. Hence
\begin{eqnarray*}
T(x)&=&T(\frac{x+y_0}{2}+\frac{x-y_0}{2})=T(\frac{x+y_0}{2})+T(\frac{x-y_0}{2})\\
&=&T(\frac{x}{2})+T(\frac{y_0}{2})+T(\frac{x}{2})+T(\frac{-y_0}{2})=2T(\frac{x}{2})+2T(\frac{y_0}{2}),
\end{eqnarray*}
\begin{eqnarray*}
T(y_0)&=&T(\frac{y_0+x}{2}+\frac{y_0-x}{2})=T(\frac{y_0+x}{2})+T(\frac{y_0-x}{2})\\
&=&T(\frac{y_0}{2})+T(\frac{x}{2})+T(\frac{y_0}{2})+T(\frac{-x}{2})=2T(\frac{y_0}{2})+2T(\frac{x}{2}),
\end{eqnarray*}
\begin{eqnarray*}
T(2x)&=&T((x+y_0)+(x-y_0))=T(x+y_0)+T(x-y_0)\\
&=&T(x)+T(y_0)+T(x)+T(-y_0)=2T(x)+2T(y_0)=4T(x),
\end{eqnarray*}
and so $T^e(2x)=4T^e(x)$. Therefore,
\begin{eqnarray*}
\|f(2x)-f(-2x)-4f(x)-4f(-x)\|&\leq&\|f^e(2x)-4f^e(x)\|\\
&=&\|f^e(2x)-T^e(2x)\|+\|4T^e(x)-4f^e(x)\|\\
&\unlhd&\varepsilon.
\end{eqnarray*}
\end{proof}

\begin{remark}
Let the binary relation $\perp'$ is defined by
\begin{eqnarray*}
x\perp' y \Leftrightarrow (x\perp y \textrm{~or~} y\perp x)
\end{eqnarray*}
Then clearly $\perp'$ is a symmetric orthogonality in the sense of
R\" atz. If $f, g, h, k$ are even mappings, then $(2.19)$ shows
that if $(2.1)$ holds for all $x,y\in {\mathcal X}$ with $x\perp
y$, then the same holds for all $x,y\in {\mathcal X}$ with
$x\perp' y$. Now if $T$ is an orthogonally additive mapping with
respect to $\perp'$ then it is trivially an orthogonally additive
mapping with respect to $\perp$. To prove the theorem therefore,
in the case that all mappings are even, we may omit the
assumption that $\perp$ is symmetric.
\end{remark}
\begin{remark} In 1985, R\"
atz (cf. Corollary 7 of \cite{RAT}) stated that if $(Y,+)$ is
uniquely $2$-divisible (i.e. the mapping $\omega:Y\to Y,
\omega(y)=2y$ is bijective), in particular $Y$ is a vector space,
then every orthogonally additive mapping $T$ has the form $T=A+P$
with $A$ additive and $P$ quadratic.
\end{remark}
The first corollary gives us a sufficient and necessary condition
to approximate an orthogonally quadratic mapping by orthogonally
additive and orthogonally quadratic mappings.
\begin{corollary} Suppose that ${\mathcal X}$ is a real orthogonality space with a symmetric orthogonal
relation $\perp$ and ${\mathcal Y}$ is a Banach space. Let
$Q:{\mathcal X}\to{\mathcal Y}$ be an orthogonally quadratic
mapping. Then a necessary and sufficient condition for the
existence of an additive mapping $A$ and an quadratic mapping $P$
with
\begin{eqnarray*}
\|Q(x)-A(x)-P(x)\|\unlhd~\varepsilon,
\end{eqnarray*}
is that
\begin{eqnarray*}
\|Q(2x)-4Q(x)\|\unlhd~\varepsilon.
\end{eqnarray*}
\end{corollary}
\begin{proof}
Set $f=g=Q$ and $h=k=2Q$ in Theorem 2.1. Then, by remark 2.3,
there exist an additive mapping $A$ and an quadratic mapping $P$
such that
\begin{eqnarray*}
\|Q(x)-A(x)-P(x)\|\unlhd~\varepsilon.
\end{eqnarray*}

Conversely, if there exists the orthogonally additive mapping
$T=A+P$ such that $\|Q(x)-T(x)\|\unlhd~\varepsilon$, then the
computations in the Step (VIII) of Theorem 2.1 gives rise
\begin{eqnarray*}
\|Q(2x)-4Q(x)\|\unlhd~\varepsilon.
\end{eqnarray*}
Note that $Q$ is orthogonally quadratic and so is clearly even,
i.e. $Q^e=Q$.
\end{proof}
The second corollary gives a solution of the stability of
Pexiderized Cauchy equation ( see also \cite{MOS1}).
\begin{corollary}
Suppose that ${\mathcal X}$ is a real orthogonality space with a
symmetric orthogonal relation $\perp$ and ${\mathcal Y}$ is a
Banach space. Let the mappings $f, h, k:{\mathcal X}\to {\mathcal
Y}$ satisfy the following inequality
\begin{eqnarray*}
\|f(x+y)-h(x)-k(y)\|\leq\varepsilon
\end{eqnarray*}
for all $x, y\in {\mathcal X}$ with $x\perp y$. Then there exists
a unique orthogonally additive mapping $T:{\mathcal X}\to
{\mathcal Y}$ such that
\begin{eqnarray*}
\|f(x)-f(0)-T(x)\|\leq 32\varepsilon
\end{eqnarray*}
\begin{eqnarray*}
\|h(x)+k(x)-h(0)-k(0)-2T(x)\|\leq 16\varepsilon
\end{eqnarray*}
for all $x\in {\mathcal X}$.
\end{corollary}
\begin{proof}
The proof of Step (IV) of Theorem 2.1 states that the condition
\begin{eqnarray*}
\|f(2x)-f(-2x)-4f(x)-4f(-x)\|\leq\unlhd~\varepsilon
\end{eqnarray*}
holds if and only if so does
\begin{eqnarray*}
\|g(2x)-g(-2x)-4g(x)-4g(-x)\|\leq\unlhd~\varepsilon.
\end{eqnarray*}
Hence we may let $G=0$ in Theorem 2.2. Then $R'=S'=0$ and the
constructions in $(2.28)$ and $(2.32)$ of the proof of the
theorem allow us to have $\|F^e(x)-S(x)\|\leq 14\varepsilon$ and
$\|L^e(x)-S(x)\|\leq 16\varepsilon$. Then
\begin{eqnarray*}
\|f(x)-f(0)-T(x)\|\leq \|F^o(x)-R(x)\|+\|F^e(x)-S(x)\|\leq
18\varepsilon+14\varepsilon=32\varepsilon,
\end{eqnarray*}
and
\begin{eqnarray*}
&&\|h(x)+k(x)-h(0)-k(0)-2T(x)\|\leq 2\|L^o(x)-R(x)\|+2\|L^e(x)-S(x)\|\\
&\leq& 40\varepsilon+32\varepsilon\\
&=&72\varepsilon.
\end{eqnarray*}
for all $x\in {\mathcal X}$.
\end{proof}
The third corollary concerns the case that ${\mathcal X}$ is
assumed to be an ordinary inner product space.
\begin{corollary} Suppose that ${\mathcal H}$ is a real inner product space of dimension greater than or equal 3 and ${\mathcal Y}$ is a Banach
space. Let the mappings $f, g, h, k:{\mathcal H}\to {\mathcal Y}$
satisfy the following inequalities
\begin{eqnarray*}
\|f(x+y)+g(x-y)-h(x)-k(y)\|\leq\varepsilon,
\end{eqnarray*}
and
\begin{eqnarray*}
\|f(2x)-f(-2x)-4f(x)-4f(-x)\|\leq\varepsilon,
\end{eqnarray*}
for all $x, y\in {\mathcal H}$ with $x\perp y$. Then there exist
orthogonally additive mappings $T, T', T'':{\mathcal X}\to
{\mathcal Y}$ such that
\begin{eqnarray*}
\|f(x)-f(0)-T(x)\|\leq \frac{140}{3}\varepsilon,
\end{eqnarray*}
\begin{eqnarray*}
\|g(x)-g(0)-T'(x)\|\leq\frac{98}{3}\varepsilon,
\end{eqnarray*}
\begin{eqnarray*}
\|h(x)+k(x)-h(0)-k(0)-T''(x)\|\leq\frac{256}{3}\varepsilon,
\end{eqnarray*}
for all $x\in {\mathcal H}$.
\end{corollary}

\end{document}